\documentclass[12pt]{article}
\usepackage{amsmath,amssymb,amscd}
\title{Publication history of von Staudt's Geometrie der Lage}
\author{Robin Hartshorne}
\date{}

\begin{document}

\maketitle

\begin{abstract}
From a census of forty copies, we can distinguish three different editions
of von Staudt's Geometrie der Lage: the first of 1847 and two undated ones
from the 1870's.
\end{abstract}

Karl Georg Christian von~Staudt's Geometrie der Lage was first published in 1847 by Bauer und Raspe in N\"urnberg, followed by three Hefts of Beitr\"age zur Geometrie der Lage in 1856, 1857, and 1860.  Raymond Clare Archibald at Brown University has already noted \cite{1} that another undated edition of the Geometrie appeared probably in the late 1870's.  But the situation is more complicated than Archibald suspected.  To clarify it, I have undertaken a census of forty copies, randomly determined by those libraries or individuals willing to return a questionnaire.

We can distinguish three editions, the first published in 1847.  The second edition is a hybrid, where pages~1--192 are identical to the first edition (these appear to be sheets remaining from the first printing) while pages~193--216 have been reset in a slightly different typeface and on a different quality of paper.  The third edition (described by Archibald) is entirely reset.  The page numbering and the text on each page are almost the same as the first edition, but the text block is slightly narrower and taller ($93 \times 165$mm compared to $96 \times 161$mm).  The three editions can be distinguished by the presence of the following errata.
\[
\begin{tabular}{|l|c|c|c|}
\hline 
\qquad\qquad\quad erratum &1st edition &2nd edition &3rd edition \\ \hline
p.~11, line~13: bdsteht & & &X \\
p.~13, last line: Ebeuen &X &X & \\
p.~103, 1st line: Din Drcieck &X &X & \\
p.~196, misnumbered 197 &X & & \\
p.~196, line~4: velcher & &X & \\
p.~203, line~13: in der Fbene PQR & &X & \\ \hline
\end{tabular}
\]

The front matter (title page, Vorwort, Inhalt) of the three editions is also different.

The first edition has on the title page Verlag von Bauer und Raspe (Julius Merz) 1847.  The reverse title page has Erlangen, gedruckt bei J.~J.~Barfus (but in 2 out of 10 copies was blank).

The hybrid second edition has variable front matter.  Three out of seven copies have the original front matter (as above).  One has the original title page, but with a paste-over Friedr.~Korn'sche Buchhandlung covering Bauer und Raspe.  Three copies have reset front matter, with Verlag der Friedr.~Korn'sche Buchhandlung (no date) on title page, and on page~vi, the line break for section~20 is Curven~II.\ ord-/nung.

The third edition has new front matter, with Verlag der Fr.~Korn'sche Buchhandllung (no date) on title page; reverse title page Druck von E.~Th.~Jacob in Erlangen; and on page~vi, the section sign (double s) for section~13 is in a different typeface.

One can also note a difference in the quality of paper.  The first edition (and first part of 2nd edition) is on whitish paper with a tendency to foxing, while the new part of the 2nd edition and the 3rd edition have a more brittle paper, with a tendency toward browning at the edges.

Of the copies examined, 10 were first editions, 7 were second editions, and 23 were third editions.  Note that having the title page of the first edition does not guarantee that it is a first edition!

I have not made a systematic attempt to study editions of the Beitr\"age, but I have seen copies saying Bauer und Raspe, with dates 1856, 1857, and 1860, presumably first editions, and other copies saying Fr.~Korn'sche Buchhandlung and no date, presumably contemporary with the 3rd edition of the Geometrie.

It is not easy to date the second and third editions of the Geometrie.  Looking at other books by the same publishers, one finds that the firm of Bauer und Raspe was active in N\"urnberg from the 1830's to the early 1870's, while Korn was active from mid-1870's to the 1930's.  In particular, W\"ockel's Geometrie der Alten \cite{9} (for which we find ads in some copies of the Geometrie) was published by Bauer und Raspe from 1839 up through the 9th edition of 1871, but from the 10th edition of 1874 on was published by Korn.  In particular, the 11th edition of 1876 has the imprint Verlag der Friedr.~Korn'sche Buchhandlung, suggesting that it is contemporary with the 2nd edition of von~Staudt described above.

One copy of the 2nd edition has a handwritten date 1878 on the title page, presumably a library acquisition date.  One copy of the 3rd edition has a previous owner's date Edinburgh 1881.  A copy of the third edition has an ad for Klingenfeld's Darstellende Geometrie in 3 volumes.  Since the third volume was first published in 1876, this places our 3rd edition after 1876.  Thus we may not be too far off if we date the 2nd edition between 1874 and 1877, and the third edition between 1877 and 1880.

From this data, one can speculate that Friedrich Korn, when taking over the business from Bauer und Raspe, found that there was a renewed demand for von~Staudt's books.  There remained sheets from the first edition, but only for pages~1--192, so he reprinted pages~193--216 and the front matter and sold more copies.  When these did not suffice, he reset the entire book, preserving the text and the page numbering and thus created our third edition.  Since the author had died in 1867, he could not ask him for a revised edition.

One can still ask, why should there be renewed interest in von~Staudt's book, thirty years after it was first published?  It seems clear that von~Staudt's work was not understood or appreciated when it first appeared.  Max Noether, in his memorial on von~Staudt \cite{5}, says that at the time of publication of the third part of the Beitr\"age in 1860, which coincided with the 25th anniversary of von~Staudt's professorship in Erlangen, his colleagues could not appreciate the value of his pathbreaking work, but respected the professor who continued his own rigorous solitary research.  Theodor Reye, whose lectures on von~Staudt's approach to projective geometry were first published in 1866 \cite{6}, says in his preface that the austere language, the extreme abstractness of presentation, and the lack of diagrams have hindered the well-deserved recognition of von~Staudt's work.  Perhaps Reye's lectures began to reawaken interest in von~Staudt.  But it seems to me that it was probably Felix Klein, with his interest in the foundations of geometry and the so-called non-Euclidean geometries, who focussed attention again on von~Staudt.  Klein in 1873 \cite{2} claimed there was a gap in von~Staudt's proof of a key result (what we now call the Fundamental Theorem of Projective Geometry), which could only be filled by an axion of continuity.  Klein's article drew responses from Cantor, L\"uroth, and Zeuthen, which Klein describes in a subsequent article of 1874 \cite{3}.  This discussion of the fundamental theorem also helped crystallize concepts of continuity, which had only been handled in a confused manner earlier.  Reye, in the preface to the second edition of his lectures in 1877, mentions Klein's objections and says that because of these he has substituted a new proof of the fundamental theorem due to Thomae.

Noether \cite[p.~116]{5} says that even von~Staudt himself did not realize the implications of his work, that it was possible to construct the metric from purely projective-geometric data.  It remained for later generations to appreciate the impact of von~Staudt's work on the foundations of projective geometry.

\bigskip
\noindent
{\bf Acknowledgements.} My thanks to all those people who patiently filled out questionnaires, answered e-mails, and sent copies of selected pages.  And special thanks to Ann Jensen of the Astronomy/Mathematics library at Berkeley, who helped diffuse the questionnaire and collect responses through her contacts in the world of librarians.


\begin{thebibliography}{99}

\bibitem{1} Archibald, R.~C., {\em Note on editions of von~Staudt's Geometrie der Lage}, Bull.~A.M.S.\ {\bf 25} (1919), 132--134.

\bibitem{2} Klein, F., {\em Ueber die sogenannte Nicht--Euklidische Geometrie (Zweiter Aufsatz)}, Math.\ Ann.\ {\bf 6} (1873), 112--145.

\bibitem{3} Klein, F., {\em Nachtrag zu dem ``Zweiten Aufsatz \"uber Nicht--Euklidische Geometrie''}, Math.\ Ann.\ {\bf 7} (1874), 531--537.

\bibitem{4} Klingenfeld, F.~A., ``Lehrbuch der Darstellenden Geometrie'', Band III, N\"urnberg, Verlag der Friedr.\ Korn'sche Buchhandlung, 1876.

\bibitem{5} Noether, M., {\em Zur Erinnerung an Karl Georg Christian von~Staudt}, Jahresbericht der Deutschen Mathematiker--Vereinigung {\bf 32} (1923), 97--119 (first published in 1901).

\bibitem{6} Reye, Th., ``Die Geometrie der Lage'', Vortr\"age (2nd ed.), Hannover, Carl R\"umpler, 1877 (first published in 1866).

\bibitem{7} Staudt, G.~K.~C.\ von, ``Geometrie der Lage'', N\"urnberg, Verlag von Bauer und Raspe (Julius Merz), 1847.

\bibitem{8} Staudt, G.~K.~C.\ von, ``Beitr\"age zur Geometrie der Lage'', N\"urnberg, Verlag von Bauer und Raspe (Julius Merz), in three parts, dated 1856, 1857, 1860, with separate title pages.

\bibitem{9} W\"ockel, L., ``Geometrie der Alten in einer Sammlung von 850 Aufgaben$\dots$ neu bearbeitet und verbessert von Th.~E.~Schroeder, elfte Auflage'', N\"urnberg, Verlag der Friedr.\ Korn'sche Buchhandlung, 1876.

\end{thebibliography}
\end{document}